\documentclass[12pt,a4paper]{amsart}
\usepackage[T1]{fontenc}
\usepackage[utf8]{inputenc}
\usepackage{lmodern}
\usepackage[margin=29mm]{geometry}
\usepackage{amsmath,amssymb,amsthm,mathtools}
\usepackage{enumitem}
\usepackage{microtype}
\usepackage[hidelinks,pdfusetitle]{hyperref}
\setlist[enumerate]{leftmargin=*,itemsep=3pt}
\newtheorem{theorem}{Theorem}[section]
\newtheorem{proposition}[theorem]{Proposition}
\newtheorem{lemma}[theorem]{Lemma}
\newtheorem{corollary}[theorem]{Corollary}
\theoremstyle{definition}
\newtheorem{definition}[theorem]{Definition}
\newtheorem{example}[theorem]{Example}
\newtheorem{remark}[theorem]{Remark}
\newcommand{\Grp}{\mathbf{Grp}}
\newcommand{\Sub}{\operatorname{Sub}}
\newcommand{\im}{\operatorname{Im}}
\newcommand{\Ker}{\operatorname{Ker}}
\newcommand{\id}{\operatorname{id}}
\newcommand{\op}{\mathrm{op}}
\newcommand{\two}{\mathbf{2}}
\newcommand{\Norm}{\operatorname{Norm}}
\newcommand{\NGrp}{\mathbf{NGrp}}
\newcommand{\NF}{\mathbf{NF}}
\newcommand{\MCL}{\mathbf{MCL}}
\newcommand{\CLm}{\mathbf{CLat}_{\wedge,\vee_{\mathrm{fin}}}^{\ne 1}}
\newcommand{\CLj}{\mathbf{CLat}_{\vee,\wedge_{\mathrm{fin}}}^{\ne 1}}
\newcommand{\Sup}{\mathbf{Sup}}
\newcommand{\eps}{\varepsilon}
\newcommand{\Z}{\mathbb Z}
\newcommand{\C}{\mathcal C}
\newcommand{\A}{\mathcal A}
\newcommand{\B}{\mathcal B}

\DeclareMathOperator{\Hom}{Hom}
\title[On proximity subgroups of groups]{On proximity subgroups of groups}
\author{Kishan Dayaram}
\author{Zurab Janelidze}
\author{Ronald Mwita Kerario}
\address[Kishan Dayaram]{Department of Mathematical Sciences, Stellenbosch University,
Private Bag X1, Matieland 7602, South Africa}
\address[Ronald Mwita Kerario]{Department of Mathematical Sciences, Stellenbosch University,
Private Bag X1, Matieland 7602, South Africa\\
African Institute for Mathematical Sciences (AIMS),
6 Melrose Road, Muizenberg, Cape Town 7945, South Africa}
\address[Zurab Janelidze]{Department of Mathematical Sciences, Stellenbosch University\\
National Institute of Theoretical and Computational Science (NITheCS)\\
Private Bag X1, Matieland 7602, South Africa}
\date{Draft --- 7 November 2026}
\subjclass[2020]{Primary 18D30; Secondary 18A32, 20J15, 06B23, 06F05}
\keywords{Noetherian form, lattice-valued norm, normed group, multiplicative complete lattice, quantale, full embedding, subgroup form, optimal noetherian form}
\thanks{Research reported in this paper as well as the presentation are each human-led and AI-assisted (ChatGPT-6 Astra). Only the second author utilized assistance of AI. All authors are responsible for the final paper draft.}
\begin{document}
\begin{abstract}
By a \emph{proximity subgroup} of a group $G$ we mean a norm $\nu$ on $G$ valued in a complete lattice $L$ equipped with a binary sum (an additive counterpart of a multiplicative structure). The idea is that the smaller the value $\nu(y)$ is, the closer $y\in G$ is to the subgroup of $G$ defined by the equation $\nu(x)=\bot$, where $\bot$ is the bottom element of $L$. We test whether this notion is a suitable substitute for the usual notion of a subgroup (recovered by taking $L=\{\bot,\top\}$ with join as the sum), by showing that the form of proximity subgroups is noetherian. Intuitively, this means that proximity subgroups are as well suited for formulating and proving isomorphism theorems and homological diagram lemmas for groups as the usual subgroups are. We study several questions in relation to this, including what makes the noetherian form of usual subgroups special among all of these noetherian forms, and in fact, all noetherian forms of the category of groups.
\end{abstract}

\maketitle

\section*{Introduction}

A noetherian form over a category is an additional structure on the category that includes a proper factorization system and enables one to formulate and prove isomorphism theorems and homological diagram lemmas 
\cite{GJ,DGJRV}.
One of the motivating examples is the \emph{form of subgroups} over the category of groups, where the isomorphism theorems and homological diagram lemmas are the usual ones. An abstract form replaces subgroups of a group $G$ with \emph{formal subgroups}, these are elements of the fibre $F^{-1}(G)$, where $F$ is a functor (the ``form'') with codomain the category of groups. In this paper we address the following question: are there other reasonable formal subgroups that are not the same as the usual subgroups, but which yield the same (or similar) isomorphism theorems and homological diagram lemmas for groups? We show that every non-trivial complete lattice $L$ equipped with a multiplicative structure (which we write in additive notation) provides such notion of a formal subgroup.

We call these formal subgroups \emph{proximity subgroups}. More precisely, an $L$-valued proximity subgroup of $G$ is a norm $\nu:G\to L$ satisfying
\[
 \nu(0)=\bot,\qquad \nu(-g)=\nu(g),\qquad
 \nu(g+h)\leq\nu(g)+\nu(h),
\]
where the coefficient sum is associative, commutative and monotone, with identity $\bot$. The norm terminology agrees with that of Clementino, Hofmann and Tholen \cite{CHT} in the quantale case. The bottom-level set of $\nu$ is an ordinary subgroup; a proximity subgroup retains the additional information supplied by its other values. This includes ordinary extended real-valued norms, as well as the non-Archimedean case in which the coefficient sum is lattice join. Note that even though we write groups additively, they are not assumed abelian.

The resulting form of proximity subgroups is fibre-complete and noetherian (in particular, it is a topological functor \cite{Brummer}). Moreover, the construction is fully faithful: morphisms between these forms are precisely those induced by maps of coefficient lattices preserving arbitrary meets, finite joins and the coefficient sum. Here morphisms of noetherian forms preserve finite fibre operations and both direct and inverse image operations, without being required to preserve arbitrary fibre meets or joins. Since every complete lattice admits join as its sum, this also yields a full embedding of a category of non-trivial complete lattices into the category of noetherian forms over groups.

The usual subgroup form has a distinguished position not only among all proximity subgroup forms, but even among all noetherian forms over the category of groups. It was established in \cite{JV} that the form of subgroups is optimal (it is isomorphic to a subform of every other noetherian form over the same base category); we prove that, up to a unique isomorphism, the form of subgroups is in fact the only optimal noetherian form over groups.

\tableofcontents

\section{Noetherian forms and their morphisms}\label{sec:forms}
We use the functorial definition of a noetherian form, together with its equivalent presentation by fibre lattices and image operations; see \cite{GJ,JV}.

\begin{definition}\label{def:form}
A \emph{form over a category $\C$} is a faithful amnestic functor
\[
                       F:\A\longrightarrow\C.
\]
Here amnestic means that an isomorphism mapped to an identity is itself an identity. We assume that every fibre is small. For $X\in\C$, the fibre $F^{-1}(X)$ is then a poset; its objects are called \emph{clusters in $X$}. We write $S\leq_f T$ if there is a morphism $S\to T$ in $\A$ mapped by $F$ to $f:X\to Y$. Such a lift, when it exists, is unique.
\end{definition}

The functor $F$ is a bifibration precisely when, for every $f:X\to Y$, there are functorial maps
\[
 f_!:F^{-1}(X)\rightleftarrows F^{-1}(Y):f^*
\]
such that
\begin{equation}\label{eq:galois}
 f_!S\leq T\quad\Longleftrightarrow\quad S\leq_f T
 \quad\Longleftrightarrow\quad S\leq f^*T.
\end{equation}
These are the direct and inverse image operations. Their functoriality includes
$(gf)_!=g_!f_!$, $(gf)^*=f^*g^*$, and identity laws. Suppose also that each fibre is a bounded lattice, with bounds $0_X,1_X$. Put
\[
 \Ker_F f=f^*(0_Y),\qquad \im_F f=f_!(1_X).
\]
A cluster is \emph{normal} if it is a kernel, and \emph{conormal} if it is an image.

For a conormal $S$ in $X$, an \emph{associated embedding} is a morphism
$i_S:S_0\to X$ with image $S$, universal among morphisms $u$ into $X$ satisfying $\im_F u\leq S$. Thus
\begin{equation}\label{eq:embedding}
 \im_F u\leq S\quad\Longleftrightarrow\quad
 u=i_Sv\ \text{for a unique }v.
\end{equation}
Dually, an \emph{associated projection} of a normal cluster $N$ in $X$ is a morphism $p_N:X\to Q_N$ with kernel $N$ such that
\begin{equation}\label{eq:projection}
 N\leq\Ker_F u\quad\Longleftrightarrow\quad
 u=vp_N\ \text{for a unique }v.
\end{equation}
These universal properties imply that associated embeddings are monomorphisms and associated projections are epimorphisms.

\begin{definition}\label{def:nf}
A form $F:\A\to\C$ is \emph{noetherian} if it is a bifibration with bounded lattice fibres and satisfies the following conditions.
\begin{enumerate}[label=\textup{(N\arabic*)}]
\item For every $f:X\to Y$ and clusters $S$ in $X$, $T$ in $Y$,
\begin{equation}\label{eq:nfidentities}
 f^*f_!S=S\vee\Ker_F f,\qquad
 f_!f^*T=T\wedge\im_F f.
\end{equation}
\item Associated embeddings and projections exist for conormal and normal clusters, respectively. Every $f$ factors as
\[
 f=i_{\im_F f}\,h\,p_{\Ker_F f},
\]
where $h$ is an isomorphism.
\item Binary joins of normal clusters are normal, and binary meets of conormal clusters are conormal.
\end{enumerate}
A noetherian form is \emph{fibre-complete} if all its fibres are complete lattices.
\end{definition}

The term \emph{noetherian} here refers to these categorical axioms, not to an ascending chain condition.

\begin{definition}\label{def:nfmorphism}
A \emph{morphism of noetherian forms} $F:\A\to\C$ and $F':\B\to\C$ is a functor $T:\A\to\B$ over $\C$ such that each fibre map preserves $0,1,\wedge,\vee$, and
\begin{equation}\label{eq:morphism}
 T_Y f_!^F=f_!^{F'}T_X,\qquad
 T_X f_F^*=f_{F'}^*T_Y.
\end{equation}
No preservation of arbitrary fibre meets or joins is required. We denote this category by $\NF(\C)$, and its full subcategory of fibre-complete objects by $\NF^c(\C)$.
\end{definition}

\begin{lemma}\label{lem:preserve}
Morphisms of noetherian forms preserve kernels, images, and their associated projections and embeddings.
\end{lemma}
\begin{proof}
Preservation of bounds and \eqref{eq:morphism} give
$T_X(\Ker_F f)=\Ker_{F'}f$ and $T_Y(\im_F f)=\im_{F'}f$.
In a noetherian form a morphism is an associated embedding if and only if its kernel is bottom; dually, it is an associated projection if and only if its image is top. These criteria follow from the factorization axiom: a projection of bottom and an embedding of top are isomorphisms. Hence an associated embedding of $S$ remains an associated embedding of $T_XS$, and dually for projections.
\end{proof}

For optimality we also need a weaker notion. An \emph{operator of forms} is any functor over the base. It is \emph{full} when
\[
 S\leq_f T\quad\Longleftrightarrow\quad
 \alpha S\leq_f\alpha T.
\]
A full operator identifies its source with a full subcategory over the base, called a \emph{subform}. In particular its fibre maps are order embeddings. Full operators are not required to preserve fibre bounds, fibre operations, or either image operation.

\begin{proposition}\label{prop:grpclasses}
Every noetherian form over $\Grp$ has all injective homomorphisms as its associated embeddings and all surjective homomorphisms as its associated projections.
\end{proposition}
\begin{proof}
Associated embeddings are injective because they are monomorphisms. Associated projections are surjective because epimorphisms in $\Grp$ are surjective. If an injective homomorphism is factored by \textup{(N2)}, its projection factor is both injective and surjective, hence an isomorphism. Thus the original homomorphism is an associated embedding. The dual argument applies to a surjective homomorphism.
\end{proof}

The \emph{subgroup form} is the functor
\[
 \mathsf S:\Grp_2\longrightarrow\Grp,\qquad (G,H)\longmapsto G,
\]
where $H\leq G$ and an arrow $(G,H)\to(K,J)$ is a homomorphism $f:G\to K$ satisfying $f(H)\leq J$. It is a fibre-complete noetherian form, with the usual direct and inverse images of subgroups.

\section{Norms with values in multiplicative complete lattices}\label{sec:norms}

\begin{definition}\label{def:coefficient}
A \emph{multiplicative complete lattice in additive convention} is a complete lattice $L$, with least element $0$ and greatest element $\infty$, equipped with an associative, commutative, monotone operation
\[
                         +:L\times L\longrightarrow L
\]
whose identity is $0$. We do not assume that $+$ distributes over meets or joins.
\end{definition}

This is the order-dual, additively written presentation of a complete lattice with monotone commutative monoid multiplication and top element as unit. Its order dual is a multiplicative lattice in the broad sense of \cite{FFJ,Facchini}, with these additional monoid and monotonicity axioms. Some authors include distribution over arbitrary joins in the multiplicative definition; the difference in terminology is discussed in \cite{GD}. In additive convention that stronger requirement is
\begin{equation}\label{eq:infdistribution}
 a+\bigwedge_i b_i=\bigwedge_i(a+b_i).
\end{equation}
It says precisely that the order dual is an integral commutative quantale. Neither \eqref{eq:infdistribution} nor any other distributivity is used in our general construction.

Monotonicity and the unit law imply
\begin{equation}\label{eq:coefficientbounds}
                a\vee b\leq a+b,\qquad a+\infty=\infty.
\end{equation}
Two examples will be used repeatedly:
\[
 ([0,\infty],\leq,+,0),\qquad
 (L,\leq,+,0)\ \text{with }a+b:=a\vee b
\]
for an arbitrary complete lattice $L$. In the first example addition is extended by $a+\infty=\infty$. In the second, coefficient addition is idempotent. Conversely, if $a+a=a$ for every $a$, then $+=\vee$: putting $c=a\vee b$ gives $c\leq a+b\leq c+c=c$.

\begin{definition}\label{def:norm}
Let $L$ be a multiplicative complete lattice in additive convention. An \emph{$L$-valued norm} on a group $G$ is a map $\nu:G\to L$ such that
\begin{equation}\label{eq:norm}
 \nu(0)=0,\qquad \nu(-g)=\nu(g),\qquad
 \nu(g+h)\leq\nu(g)+\nu(h).
\end{equation}
The pair $(G,\nu)$ is called an \emph{$L$-normed group}.
A homomorphism $f:(G,\nu)\to(H,\mu)$ is \emph{nonexpansive}, or an \emph{$L$-normed homomorphism}, when
\begin{equation}\label{eq:nonexpansive}
                         \mu(f(g))\leq\nu(g)\qquad(g\in G).
\end{equation}
These objects and arrows form a category $\NGrp_L$, with forgetful functor
\[
                         U_L:\NGrp_L\longrightarrow\Grp.
\]
We write $\Norm_L(G)$ for its fibre. Its order is reverse pointwise order:
\begin{equation}\label{eq:reverse}
           \nu\preceq \mu\quad\Longleftrightarrow\quad
           \mu(g)\leq\nu(g)\quad\text{for every }g\in G.
\end{equation}
\end{definition}

We use \emph{norm} in the categorical sense of Clementino, Hofmann and Tholen \cite[\S5]{CHT}: symmetry under inversion is required, but separation and finiteness are not. Thus $\nu(g)=0$ need not imply $g=0$, and the value $\infty$ is allowed. When \eqref{eq:infdistribution} holds, $\mathcal V=L^{\op}$ is an integral commutative quantale, and our $L$-normed groups are precisely their $\mathcal V$-normed groups. Our nonexpansive homomorphisms correspond to $\mathcal V$-normed functors between the associated one-object categories. We retain these group-norm axioms for the more general coefficients of Definition~\ref{def:coefficient}, without assuming distributivity.

Coefficient comparisons use $\leq$, whereas fibre comparisons use $\preceq$.
We distinguish fibre operations by $\wedge_F,\vee_F$ when needed.
Ordinary extended real-valued group norms are obtained with $L=[0,\infty]$ and ordinary addition; no homogeneity or conjugation invariance is imposed. For related real-valued length functions see \cite{Polymath}. Such norms correspond to maps $S:G\to[0,1]$ with
\[
 S(0)=1,\qquad S(-g)=S(g),\qquad S(g)S(h)\leq S(g+h),
\]
by $S(g)=e^{-\nu(g)}$, with $e^{-\infty}=0$ and $-\log0=\infty$.
With $a+b=a\vee b$, condition \eqref{eq:norm} instead becomes the non-Archimedean inequality. For real coefficients this uses maximum in place of addition. The lattice-valued fuzzy-subgroup literature, beginning with \cite{Goguen,Rosenfeld}, concerns the order-dual meet inequality.

\begin{lemma}\label{lem:bifibration}
The functor $U_L$ is a bifibration with complete lattice fibres. Inverse image is
\[
                             f^*\mu=\mu\circ f.
\]
For an arbitrary map $u:G\to L$, define
\begin{equation}\label{eq:closure}
 c_G(u)(g)=
 \bigvee\{\eta(g):\eta\in\Norm_L(G),\ \eta(x)\leq u(x)\text{ for all }x\in G\}.
\end{equation}
Then $c_G(u)$ is the greatest norm pointwise below $u$, and direct image is
\begin{equation}\label{eq:direct}
 (f_!\nu)(h)=c_H(r_f\nu)(h),\qquad
 (r_f\nu)(h)=\bigwedge_{\{g:f(g)=h\}}\nu(g).
\end{equation}
The empty meet in \eqref{eq:direct} is $\infty$.
\end{lemma}
\begin{proof}
Faithfulness and amnesticity are immediate. Arbitrary pointwise suprema of norms are norms. Indeed,
\[
 \bigvee_i\nu_i(g+h)
 \leq\bigvee_i\bigl(\nu_i(g)+\nu_i(h)\bigr)
 \leq\left(\bigvee_i\nu_i(g)\right)+
      \left(\bigvee_i\nu_i(h)\right).
\]
The empty supremum is the zero function. These suprema are the fibre meets:
\[
                  \left(\bigwedge_{F,i}\nu_i\right)(g)
                          =\bigvee_i\nu_i(g).
\]
Thus all meets exist and the fibres are complete. Their joins are
\[
                         \bigvee_{F,i}\nu_i
                         =c_G\!\left(\bigwedge_i\nu_i\right).
\]
The zero norm lies pointwise below every $u$, and the preceding observation proves the assertion about $c_G$.

Composition with a homomorphism preserves norms. For $\mu\in\Norm_L(H)$,
\[
 \begin{aligned}
 f_!\nu\preceq \mu
 &\Longleftrightarrow \mu\leq c_H(r_f\nu)\quad\text{pointwise}\\
 &\Longleftrightarrow \mu\leq r_f\nu\\
 &\Longleftrightarrow \mu\circ f\leq\nu\\
 &\Longleftrightarrow \nu\preceq f^*\mu.
 \end{aligned}
\]
This is the adjunction describing lifts of $f$. Inverse images are functorial; uniqueness of left adjoints gives functoriality of direct images.
\end{proof}

For a subgroup $A\leq G$, set
\begin{equation}\label{eq:characteristic}
 \kappa_A(g)=
 \begin{cases}
 0,&g\in A,\\
 \infty,&g\notin A.
 \end{cases}
\end{equation}
These are norms by \eqref{eq:coefficientbounds}. The fibre bottom and top are $\kappa_{\{0\}}$ and $\kappa_G$, respectively.

\begin{lemma}\label{lem:intervals}
Let $p:G\twoheadrightarrow H$ be surjective with kernel $N$, and let $i:A\hookrightarrow G$ be an inclusion.
\begin{enumerate}[label=\textup{(\roman*)}]
\item Inverse image along $p$ identifies $\Norm_L(H)$ with the upper fibre interval $\uparrow\kappa_N$. Consequently
\[
                  p_!p^*=\id,\qquad
                  p^*p_!\nu=\nu\vee_F\kappa_N.
\]
\item Direct image along $i$ extends a norm on $A$ by $\infty$ outside $A$. It identifies $\Norm_L(A)$ with the lower fibre interval $\downarrow\kappa_A$. Consequently
\[
                  i^*i_!=\id,\qquad
                  i_!i^*\nu=\nu\wedge_F\kappa_A.
\]
\end{enumerate}
\end{lemma}
\begin{proof}
The condition $\kappa_N\preceq\nu$ says exactly that $\nu$ vanishes on $N$.
For $n\in N$ the triangle inequality gives
\[
 \nu(g+n)\leq\nu(g)+\nu(n)=\nu(g),\qquad
 \nu(g)\leq\nu(g+n)+\nu(-n)=\nu(g+n).
\]
Hence $\nu$ is constant on the cosets of $N$, and descends uniquely to a norm on $H$. Conversely every pullback vanishes on $N$. The reflector onto the upper interval is fibre join with $\kappa_N$, proving (i).

Extension by $\infty$ is a norm: if both summands lie in $A$, use the inequality in $A$; otherwise its right-hand side is $\infty$. A norm lies below $\kappa_A$ in fibre order precisely when it is $\infty$ outside $A$. Restriction and extension give the stated interval isomorphism, proving (ii).
\end{proof}

\begin{theorem}\label{thm:nf}
For every nontrivial multiplicative complete lattice $L$ in additive convention, the functor $U_L$ is a fibre-complete noetherian form. For every homomorphism $f:G\to H$,
\begin{equation}\label{eq:kernelsimages}
                  \Ker_{U_L}f=\kappa_{\ker f},\qquad
                  \im_{U_L}f=\kappa_{f(G)}.
\end{equation}
Its associated embeddings and projections are precisely the injective and surjective group homomorphisms, respectively.
\end{theorem}
\begin{proof}
Factor $f$ as $G\xrightarrow{p}f(G)\xrightarrow{i}H$.
Lemma~\ref{lem:intervals} gives
\[
 \begin{aligned}
 f^*f_!\nu&=p^*i^*i_!p_!\nu
            =p^*p_!\nu=\nu\vee_F\kappa_{\ker f},\\
 f_!f^*\mu&=i_!p_!p^*i^*\mu
            =i_!i^*\mu=\mu\wedge_F\kappa_{f(G)}.
 \end{aligned}
\]
The fibre bounds give \eqref{eq:kernelsimages}, proving \textup{(N1)}.

Since $0\ne\infty$ in $L$,
\begin{equation}\label{eq:reflectsubgroups}
                         \kappa_A\preceq\kappa_B
                         \quad\Longleftrightarrow\quad A\leq B.
\end{equation}
Thus $A\hookrightarrow G$ has image $\kappa_A$ and its ordinary subgroup universal property is exactly \eqref{eq:embedding}. For $N\trianglelefteq G$, the quotient $G\twoheadrightarrow G/N$ has kernel $\kappa_N$ and satisfies \eqref{eq:projection}. The usual group factorization establishes \textup{(N2)}. The assertion about all associated maps also follows from Proposition~\ref{prop:grpclasses}.

Finally, for normal subgroups $N,M$ and arbitrary subgroups $A,B$,
\begin{equation}\label{eq:normaloperations}
 \kappa_N\vee_F\kappa_M=\kappa_{N+M},\qquad
 \kappa_A\wedge_F\kappa_B=\kappa_{A\cap B}.
\end{equation}
The second equality is pointwise. For the first, the zero set of any norm is a subgroup. A fibre upper bound of $\kappa_N,\kappa_M$ therefore vanishes on the subgroup $N+M$ they generate, and hence lies above $\kappa_{N+M}$ in fibre order. This proves \textup{(N3)}.
\end{proof}

\begin{corollary}\label{cor:normalconormal}
The normal clusters of $U_L$ over $G$ are exactly $\kappa_N$ with $N\trianglelefteq G$, and the conormal clusters are exactly $\kappa_A$ with $A\leq G$. Every normal cluster is conormal.
\end{corollary}
\begin{proof}
Use the quotient $G\twoheadrightarrow G/N$, the inclusion $A\hookrightarrow G$, and \eqref{eq:kernelsimages}.
\end{proof}

\begin{proposition}\label{prop:characteristic}
The assignments $A\mapsto\kappa_A$ define a full operator $\mathsf S\to U_L$ that is also a morphism of noetherian forms. It preserves all fibre meets and joins.
\end{proposition}
\begin{proof}
For subgroups $A_i\leq G$,
\[
 \bigwedge_{F,i}\kappa_{A_i}=\kappa_{\bigcap_iA_i},\qquad
 \bigvee_{F,i}\kappa_{A_i}=\kappa_{\langle\bigcup_iA_i\rangle}.
\]
The meet formula is pointwise. The join formula follows because a norm vanishing on all $A_i$ vanishes on the subgroup they generate. These arguments include empty families. Moreover,
\[
 f^*\kappa_B=\kappa_{f^{-1}(B)},\qquad
 f_!\kappa_A=\kappa_{f(A)}.
\]
The first equality is immediate. The second follows from the adjunction and the fact that the zero set of a norm is a subgroup.
Finally, $\kappa_A\leq_f\kappa_B$ is equivalent to $f(A)\leq B$, which proves fullness.
\end{proof}

Write $\two=\{0,\infty\}$ with addition given by maximum. Then $U_{\two}\cong\mathsf S$.
For the cyclic group $C_2=\{0,t\}$, every $a\in L$ occurs uniquely as $\nu(t)$; hence
\begin{equation}\label{eq:C2}
                              \Norm_L(C_2)\cong L^{\op}.
\end{equation}
If $0<a<\infty$, this norm on $C_2$ is neither normal nor conormal.
Nontriviality is essential: for the singleton coefficient lattice, every homomorphism has bottom kernel and top image, and the factorization axiom would force every homomorphism to be an isomorphism.

\begin{proposition}\label{prop:quantaledirect}
If \eqref{eq:infdistribution} holds, then direct image requires no closure:
\[
                    (f_!\nu)(h)=\bigwedge_{f(g)=h}\nu(g).
\]
In particular, for ordinary extended norms,
\[
                    (f_!\nu)(h)=\inf_{f(g)=h}\nu(g).
\]
\end{proposition}
\begin{proof}
Put $u=r_f\nu$. It has $u(0)=0$ and $u(-h)=u(h)$. Distributing addition over both infima gives
\[
 \begin{aligned}
 u(h)+u(k)
 &=\bigwedge_{\substack{f(g)=h\\f(g')=k}}
       \bigl(\nu(g)+\nu(g')\bigr)\\
 &\geq\bigwedge_{\substack{f(g)=h\\f(g')=k}}\nu(g+g')
 \geq\bigwedge_{f(z)=h+k}\nu(z)=u(h+k).
 \end{aligned}
\]
Thus $u$ is already a norm. Extended real addition satisfies \eqref{eq:infdistribution}.
\end{proof}

\begin{example}\label{ex:closure}
The closure in \eqref{eq:direct} cannot be omitted in general. Let $L=M_3$ be the five-element lattice with distinct atoms $a,b,c$, and use coefficient addition $+=\vee$.
In $G=(C_2)^4$, with basis $e_1,e_2,e_3,e_4$, put
\[
 A_a=\langle e_1\rangle,\qquad
 A_b=\langle e_2,e_3\rangle,\qquad
 A_c=\langle e_4\rangle.
\]
Define $\nu(0)=0$, give the nonzero elements of $A_a,A_b,A_c$ the values $a,b,c$, respectively, and give every other element value $\infty$. This is a norm: each set with a fixed atomic upper bound is a subgroup, and the join of distinct atoms is $\infty$.
Let $q:G\to(C_2)^2$ send $(u,v,w,z)$ to $(u+v,w+z)$, and write $x=(1,0)$, $y=(0,1)$. Its raw fibre infimum satisfies
\[
                 (r_q\nu)(x)=0,\qquad
                 (r_q\nu)(y)=0,\qquad
                 (r_q\nu)(x+y)=b.
\]
It fails the triangle inequality. Its norm closure is the zero function, which is the actual direct image.
\end{example}

\begin{proposition}\label{prop:cuts}
For $a\in L$, the sublevel set
\[
                              \{g:\nu(g)\leq a\}
\]
is a subgroup for every $L$-valued norm on every group if and only if $a+a=a$. Consequently all sublevel sets of all such norms are subgroups if and only if coefficient addition is lattice join.
\end{proposition}
\begin{proof}
If $a+a=a$, the identity, inverse and sum axioms establish the subgroup condition.
Conversely, on $C_2\oplus C_2=\{0,x,y,x+y\}$ define
\[
                    \nu(0)=0,\quad
                    \nu(x)=\nu(y)=a,\quad
                    \nu(x+y)=a+a.
\]
This is a norm; the remaining inequalities follow from monotonicity and $0\leq a$. If $a+a\ne a$, its $a$-sublevel set contains $x$ and $y$ but not $x+y$. The final assertion follows from the idempotence characterization after \eqref{eq:coefficientbounds}.
\end{proof}

For example, $\nu(n)=|n|$ on $\Z$ has the sublevel set $\{-1,0,1\}$ at $a=1$, which is not a subgroup. This distinction does not affect Theorem~\ref{thm:nf}, Corollary~\ref{cor:normalconormal}, or the characteristic embedding.

\section{A full embedding of multiplicative complete lattices}\label{sec:full}

Let $\MCL$ denote the category of nontrivial multiplicative complete lattices in additive convention. A morphism $h:L\to M$ preserves arbitrary meets, finite joins, and addition:
\begin{equation}\label{eq:coefficientmap}
 \begin{gathered}
 h\left(\bigwedge_i a_i\right)=\bigwedge_i h(a_i),\qquad
 h(a\vee b)=h(a)\vee h(b),\\
 h(0)=0,\qquad h(a+b)=h(a)+h(b).
 \end{gathered}
\end{equation}
The empty meet gives $h(\infty)=\infty$. Thus bounds are preserved. In the order-dual multiplicative convention, these are the maps preserving arbitrary joins, finite meets, and multiplication. The choice of morphisms will be recovered from the noetherian forms themselves.

\begin{proposition}\label{prop:coefficient}
For every morphism $h:L\to M$ in $\MCL$, postcomposition defines a morphism of noetherian forms
\[
 P(h):U_L\longrightarrow U_M,\qquad
 P(h)_G(\nu)=h\circ\nu.
\]
These assignments define a functor
\[
                            P:\MCL\longrightarrow\NF^c(\Grp).
\]
Every $P(h)_G$ also preserves arbitrary fibre joins.
\end{proposition}
\begin{proof}
Preservation of $0$ and addition, together with monotonicity, ensures that $h\nu$ is a norm. Postcomposition preserves nonexpansive homomorphisms and commutes with inverse image.

Since $h$ preserves arbitrary meets, it has a left adjoint $\alpha:M\to L$:
\[
                  \alpha(b)\leq a\quad\Longleftrightarrow\quad b\leq h(a).
\]
The map $\alpha$ preserves $0$ and is subadditive. Indeed, the unit inequalities give
\[
 b+c\leq h(\alpha b)+h(\alpha c)=h(\alpha b+\alpha c),
\]
and hence
\[
                             \alpha(b+c)\leq\alpha b+\alpha c.
\]
Consequently $\alpha \mu$ is a norm whenever $\mu$ is one.
Writing $\Phi_G(\nu)=h\nu$ and $A_G(\mu)=\alpha \mu$, the pointwise adjunction becomes
\[
 \Phi_G(\nu)\preceq \mu
 \ \Longleftrightarrow\ \mu\leq h\nu
 \ \Longleftrightarrow\ \alpha \mu\leq\nu
 \ \Longleftrightarrow\ \nu\preceq A_G(\mu).
\]
Thus $\Phi_G\dashv A_G$ in fibre order, so $\Phi_G$ preserves arbitrary fibre joins.
Finite fibre meets are pointwise finite coefficient joins, which $h$ preserves.
Both fibre bounds are preserved because $h(0)=0$ and $h(\infty)=\infty$.

Postcomposition with $\alpha$ commutes with inverse image:
\[
                         A_G f_M^*=f_L^* A_H.
\]
Uniqueness of left adjoints therefore gives
\[
                         \Phi_H f_!^L=f_!^M\Phi_G.
\]
This establishes compatibility with direct image. Identities and composition are immediate.
\end{proof}

\begin{lemma}\label{lem:pointwise}
Let $T:U_L\to U_M$ be a morphism of noetherian forms. There is a unique bounded lattice homomorphism $h:L\to M$ such that
\begin{equation}\label{eq:pointwise}
                             (T_G\nu)(g)=h(\nu(g))
\end{equation}
for every group $G$, norm $\nu$, and element $g$.
\end{lemma}
\begin{proof}
By \eqref{eq:C2}, the fibre map $T_{C_2}:L^{\op}\to M^{\op}$ determines a bounded lattice homomorphism $h:L\to M$. We prove that it determines every component.

Let $\eta:\Z\to L$ be a norm. Write the infinite dihedral group additively as
\[
 D_\infty=\Z\rtimes C_2,\qquad
 (n,\eps)+(m,\delta)=
 \bigl(n+(-1)^\eps m,\eps+\delta\bigr).
\]
Here $\eps,\delta\in\{0,1\}$ and their sum is modulo $2$. Put
\[
 r=(1,0),\qquad s=(0,1),\qquad
 \lambda(n,\eps)=\eta(n).
\]
Symmetry and the triangle inequality for $\eta$ show that $\lambda$ is a norm.
Both $s$ and $r+s$ have order $2$. Restricting to the two corresponding homomorphisms from $C_2$ gives
\[
                (T_{D_\infty}\lambda)(s)=h(0)=0,\qquad
                (T_{D_\infty}\lambda)(r+s)=h(\eta(1)).
\]
For any norm $\omega$ vanishing at $s$, the triangle inequality and $s+s=0$ imply
\[
                  \omega(r+s)\leq\omega(r),\qquad
                  \omega(r)\leq\omega(r+s).
\]
Thus $(T_{D_\infty}\lambda)(r)=h(\eta(1))$. Pulling back along
$\Z\to D_\infty$, $n\mapsto(n,0)$, yields
\[
                              (T_\Z\eta)(1)=h(\eta(1)).
\]
Finally, pull back an arbitrary $\nu$ along $u_g:\Z\to G$, $1\mapsto g$.
Naturality of inverse image gives
\[
 (T_G\nu)(g)
 =\bigl(T_\Z(u_g^*\nu)\bigr)(1)
 =h(\nu(g)).
\]
Uniqueness already follows from the fibre over $C_2$.
\end{proof}

\begin{lemma}\label{lem:addition}
The map $h$ in Lemma~\ref{lem:pointwise} preserves coefficient addition.
\end{lemma}
\begin{proof}
Let $V=C_2\oplus C_2=\{0,x,y,x+y\}$. For $a,b\in L$, list values at these four elements and consider the norms
\[
                         A=(0,a,\infty,\infty),\qquad
                         B=(0,\infty,b,\infty).
\]
Their fibre join is
\begin{equation}\label{eq:V4join}
                             A\vee_F B=(0,a,b,a+b).
\end{equation}
The displayed tuple is a norm. Its triangle inequalities follow from commutativity, associativity and \eqref{eq:coefficientbounds}; for example, the only bound involving both independent generators is equality at $x+y$, and $a\leq b+(a+b)$ follows from monotonicity.
Any norm pointwise below both $A$ and $B$ has values at most $a,b,a+b$ at $x,y,x+y$, respectively. This proves \eqref{eq:V4join}.

Apply $T_V$, which preserves this fibre join, and use \eqref{eq:pointwise}. Evaluation at $x+y$ gives
\[
                                  h(a+b)=h(a)+h(b).
\]
\end{proof}

\begin{lemma}\label{lem:arbitrarymeets}
The map $h$ in Lemma~\ref{lem:pointwise} preserves arbitrary meets.
\end{lemma}
\begin{proof}
Let $(a_i)_{i\in I}$ be a family in $L$, and identify elements of
$E_I=\bigoplus_{i\in I}C_2$ with finite subsets $J\subseteq I$, using symmetric difference as addition. Define
\[
                      \nu_I(J)=\sum_{j\in J}a_j,\qquad
                      \nu_I(\varnothing)=0.
\]
This is a norm. Passing from the concatenated sums for $J$ and $K$ to the sum for $J\mathbin{\triangle}K$ removes duplicate nonnegative summands, and hence cannot increase the value.

Let $\eps_I:E_I\to C_2$ be parity. Direct image into $C_2$ requires no correction of the raw fibre meet: every value at its nonzero element $t$, together with value $0$ at $0$, defines a norm. Thus
\begin{equation}\label{eq:parity}
  (\eps_{I!}\nu_I)(t)
  =\bigwedge_{\substack{J\subseteq I\ \mathrm{finite}\\|J|\ \mathrm{odd}}}
          \sum_{j\in J}a_j
  =\bigwedge_{i\in I}a_i.
\end{equation}
Every nonempty sum is above each of its summands, while singleton subsets give the converse bound. If $I$ is empty, both sides are $\infty$.

By Lemma~\ref{lem:addition}, $h$ preserves finite sums. Using \eqref{eq:pointwise}, commutation with direct image, and \eqref{eq:parity}, we obtain
\[
 \begin{aligned}
 h\left(\bigwedge_i a_i\right)
 &=\bigl(T_{C_2}(\eps_{I!}\nu_I)\bigr)(t)\\
 &=\bigl(\eps_{I!}(T_{E_I}\nu_I)\bigr)(t)\\
 &=\bigwedge_{\substack{J\subseteq I\ \mathrm{finite}\\|J|\ \mathrm{odd}}}
                  \sum_{j\in J}h(a_j)
  =\bigwedge_i h(a_i).
 \end{aligned}
\]
\end{proof}

\begin{theorem}\label{thm:full}
The functor
\[
                              P:\MCL\longrightarrow\NF^c(\Grp)
\]
is full and faithful. Explicitly,
\[
 \Hom_{\MCL}(L,M)
 \ \cong\
 \Hom_{\NF(\Grp)}(U_L,U_M),\qquad
 h\longmapsto(\nu\mapsto h\nu).
\]
\end{theorem}
\begin{proof}
Proposition~\ref{prop:coefficient} constructs a morphism of noetherian forms from each coefficient morphism. Lemmas~\ref{lem:pointwise}, \ref{lem:addition} and \ref{lem:arbitrarymeets} recover a unique coefficient morphism from every morphism of forms.
\end{proof}

\begin{corollary}\label{cor:latticeembedding}
Let $\CLm$ be the category of nontrivial complete lattices and maps preserving arbitrary meets and finite joins. The assignment
\[
                         L\longmapsto U_{(L,+\,=\,\vee)}
\]
is full and faithful. Equivalently, for the category $\CLj$ of nontrivial complete lattices and maps preserving arbitrary joins and finite meets, the assignment
\[
                         A\longmapsto U_{(A^{\op},+\,=\,\vee_{A^{\op}})}
\]
is full and faithful.
\end{corollary}
\begin{proof}
The idempotent objects of $\MCL$ are exactly complete lattices with addition equal to join. Preservation of finite joins already ensures preservation of this addition, so these objects form the stated full subcategory. The second assertion follows by order duality.
\end{proof}

\begin{corollary}\label{cor:range}
The forms $U_L$ and $U_M$ are isomorphic over $\Grp$ if and only if $L$ and $M$ are isomorphic as multiplicative complete lattices. These forms have arbitrarily large fibres already over $C_2$.
\end{corollary}
\begin{proof}
Full faithfulness reflects and preserves isomorphisms, and \eqref{eq:C2} identifies the fibre over $C_2$.
\end{proof}

In particular the underlying complete lattice alone does not determine the form in the general construction. On $[0,\infty]$, ordinary addition and maximum yield nonisomorphic forms: the second addition is idempotent, while the first is not, so their coefficient structures cannot be isomorphic.

\begin{remark}\label{rem:infinite}
Although morphisms of noetherian forms need only preserve finite fibre operations, every morphism between the forms $U_L$ preserves arbitrary fibre joins, by Theorem~\ref{thm:full} and Proposition~\ref{prop:coefficient}. It need not preserve arbitrary fibre meets.
For example, give $[0,\infty]$ coefficient addition $a+b=\max(a,b)$, and let $\two$ have the same addition. The map
\[
 h:[0,\infty]\longrightarrow\two,\qquad
 h(t)=
 \begin{cases}
 0,&t<\infty,\\
 \infty,&t=\infty,
 \end{cases}
\]
preserves arbitrary meets, finite joins and addition. It fails to preserve the supremum of the sequence $1,2,3,\ldots$. Since fibre meets over $C_2$ are coefficient joins, the induced morphism fails to preserve arbitrary fibre meets.
\end{remark}

\section{The unique optimal form over groups}\label{sec:optimal}

We now use full operators in the weaker sense introduced in Section~\ref{sec:forms}. Following \cite{JV}, a noetherian form is \emph{optimal} if it is isomorphic to a subform of every noetherian form over the same base inducing the same associated projection--embedding factorization system. The subgroup form is optimal by \cite{JV}. By Proposition~\ref{prop:grpclasses}, all noetherian forms over $\Grp$ induce the same factorization system.

\begin{lemma}[Rigidity of the subgroup form]\label{lem:rigidity}
Every full operator $T:\mathsf S\to\mathsf S$ over $\Grp$ is the identity.
\end{lemma}
\begin{proof}
Fullness says that, for every $f:G\to K$, $H\leq G$ and $J\leq K$,
\begin{equation}\label{eq:fulloperator}
 f(H)\leq J
       \quad\Longleftrightarrow\quad
 f(T_G(H))\leq T_K(J).
\end{equation}
In particular each $T_G$ is an order embedding. For every prime $p$, the fibre over $C_p$ is the two-element chain. Hence $T_{C_p}$ fixes both subgroups.

Write $T_\Z(\Z)=n\Z$ with $n\geq0$. For the quotient $q_p:\Z\to C_p$, equation \eqref{eq:fulloperator}, applied to $\Z$ and the trivial subgroup of $C_p$, implies
\[
                           q_p(n\Z)\nleq0.
\]
Thus $n\neq0$ and no prime divides $n$, so $n=1$.

For any $g\in G$, let $u_g:\Z\to G$ send $1$ to $g$. For every $H\leq G$,
\[
 \begin{aligned}
 g\in H
 &\Longleftrightarrow u_g(\Z)\leq H\\
 &\Longleftrightarrow u_g(T_\Z(\Z))\leq T_G(H)\\
 &\Longleftrightarrow g\in T_G(H).
 \end{aligned}
\]
Therefore $T_G(H)=H$ for every subgroup of every group. Since the subgroup form is faithful, $T$ is also the identity on arrows.
\end{proof}

\begin{theorem}\label{thm:optimal}
The subgroup form is the only optimal noetherian form over $\Grp$, up to a unique isomorphism over $\Grp$.
\end{theorem}
\begin{proof}
Let $F$ be optimal. Optimality of $F$ and of $\mathsf S$ gives full operators
\[
                       \mathsf S\xrightarrow{a}F
                           \xrightarrow{b}\mathsf S.
\]
Lemma~\ref{lem:rigidity} gives $ba=\id_{\mathsf S}$. Since a full operator between forms is injective on each fibre, the equality
$b(ab)=b$ implies $ab=\id_F$. Thus $a,b$ are inverse isomorphisms. Any two isomorphisms $F\to\mathsf S$ differ by an automorphism of $\mathsf S$, which is the identity by the same lemma.
\end{proof}

\begin{corollary}\label{cor:conormal}
For a noetherian form $F$ over $\Grp$, the following conditions are equivalent:
\begin{enumerate}[label=\textup{(\roman*)}]
\item $F$ is optimal;
\item $F$ is isomorphic to the subgroup form;
\item every cluster of $F$ is conormal.
\end{enumerate}
\end{corollary}
\begin{proof}
The first two conditions are equivalent by Theorem~\ref{thm:optimal}; the subgroup form is conormal. If every cluster of $F$ is conormal, associated embeddings identify its clusters with ordinary subgroups by Proposition~\ref{prop:grpclasses}. More explicitly,
$H\leq G$ corresponds to $\im_F(H\hookrightarrow G)$. The universal property of associated embeddings identifies order and the relation along every homomorphism with those of the subgroup form. This is an isomorphism of forms.
\end{proof}

\begin{corollary}\label{cor:optimalcoefficient}
For a nontrivial multiplicative complete lattice $L$ in additive convention, the form $U_L$ is optimal if and only if $L\cong\two$ as a coefficient structure.
\end{corollary}
\begin{proof}
If $U_L$ is optimal, then it is isomorphic to $\mathsf S$. Evaluating at $C_2$ gives $L^{\op}\cong\two$. On a two-element coefficient lattice addition is forced to be join, so the converse follows from Proposition~\ref{prop:characteristic}.
\end{proof}

\begin{remark}
Theorem~\ref{thm:optimal} is stronger than the existence of full operators in both directions: the rigidity lemma is what makes those operators inverse. It also implies that the subgroup form has no proper subform which is itself a noetherian form over groups.

Optimality does not assert initiality in $\NF(\Grp)$. The operators appearing in its definition need not preserve bounds or image operations. For instance,
\[
 D(G)=\Sub(G)\times\two,
\]
with both transports acting identically on the second coordinate, is a noetherian form. Its kernels are $(\ker f,0)$, its images are $(f(G),\infty)$, and the usual group factorizations verify the axioms. There can be no morphism $\mathsf S\to D$ in $\NF(\Grp)$: at the trivial group such a morphism would have to be a bounded lattice homomorphism from the singleton lattice to $\two$.
\end{remark}

\section{Future directions}\label{sec:future}

\medskip
\noindent\textbf{Tensor product of lattices.}
The full embedding includes both ordinary extended norms and the idempotent coefficient case. The latter retains a tensor description based solely on the subgroup form. Let $\Sup$ denote the category of complete lattices and arbitrary-join preserving maps. Its tensor product $A\otimes_{\Sup}B$ represents maps $A\times B\to C$ preserving arbitrary joins separately in each variable; see \cite{GHK}.

Suppose coefficient addition on $L$ is lattice join. The sublevel subgroups of $\nu\in\Norm_L(G)$ define a map
\[
                   H_\nu:L\longrightarrow\Sub(G),\qquad
                   H_\nu(a)=\{g:\nu(g)\leq a\},
\]
with
\[
                   H_\nu\left(\bigwedge_i a_i\right)
                         =\bigcap_iH_\nu(a_i).
\]
Conversely a map preserving arbitrary meets in this way comes from a unique norm:
\[
                          \nu(g)=\bigwedge\{a:g\in H(a)\}.
\]
Preservation of arbitrary meets ensures that the displayed infimum is an attained level; the subgroup condition yields the join triangle inequality. In the join-reversing description of the tensor product this gives the fibrewise lattice isomorphism
\begin{equation}\label{eq:tensor}
                    \Norm_L(G)\cong L^{\op}\otimes_{\Sup}\Sub(G).
\end{equation}

For a fibre-complete noetherian form $F$ over another base category, the corresponding candidate for tensoring with a complete lattice $A$ has fibres
\[
                    (A\odot F)^{-1}(X)
                           =A\otimes_{\Sup}F^{-1}(X),
\]
with direct images induced by $1_A\otimes f_!$.
A systematic study should determine the hypotheses under which these data satisfy all noetherian axioms, and the appropriate class of morphisms for which tensoring is functorial. These are separate questions: the tensor product is defined using arbitrary joins, whereas morphisms of noetherian forms need not preserve arbitrary joins. The full embedding established here does not, on its own, provide an action on the entire category of noetherian forms.

For general coefficient addition, Proposition~\ref{prop:cuts} identifies the precise obstruction to the preceding level-set description. Only idempotent thresholds are guaranteed to define subgroups. Consequently \eqref{eq:tensor} must be restricted to the idempotent case. A tensor construction for the full class of norms would have to retain coefficient addition as part of its data: the underlying complete lattice alone does not determine the form, as Corollary~\ref{cor:range} shows. This issue is independent of the closure needed for direct images; Proposition~\ref{prop:quantaledirect} removes that closure for quantale coefficients, including ordinary extended real-valued norms, without making all sublevel sets subgroups.

\medskip
\noindent\textbf{Beyond groups.}
For other group-like structures, one may impose analogous inequalities involving their operations and coefficient addition. The questions are then whether the resulting forgetful functors are noetherian and whether the coefficient structures can be reconstructed from morphisms of forms. Our reconstruction uses $C_2$, the infinite dihedral group, $C_2\oplus C_2$, and parity maps. Suitable replacements and the additional properties of the starting form must be identified for each proposed generalization. These extensions and a general tensor formulation are left for future work.

\medskip
\noindent\textbf{Arbitrary pseudometrics.}
A further direction is to replace norms by arbitrary $L$-valued pseudometrics on the underlying group. These are functions $d:G\times G\to L$ satisfying
\[
 d(x,x)=0,\qquad d(x,y)=d(y,x),\qquad
 d(x,z)\leq d(x,y)+d(y,z).
\]
No separation, finiteness or compatibility with the group operation is imposed. With nonexpansive group homomorphisms they form a category $\mathbf{PMetGrp}_L$, whose forgetful functor we denote by $V_L$. Norms identify with the left-translation-invariant pseudometrics through
\[
 d_\nu(x,y)=\nu(-x+y),\qquad \nu(g)=d_\nu(0,g).
\]
When $L=\two$, arbitrary pseudometrics are equivalence relations: $x\sim y$ precisely when $d(x,y)=0$. Thus $V_{\two}$ is the pullback of the equivalence-relation form on pointed sets along $\Grp\to\mathbf{Set}_*$, as in van Niekerk's construction \cite{vN}. It differs from $U_{\two}$, the subgroup form.

The basic construction extends beyond the invariant case. Pointwise suprema of pseudometrics are pseudometrics, so the reverse pointwise fibres are complete. Let $c_G^{\mathrm{met}}(u)$ denote the greatest pseudometric pointwise below a function $u:G\times G\to L$. Then
\[
 f^*e(x,y)=e(fx,fy),\qquad
 f_!d=c_H^{\mathrm{met}}(r_fd),\qquad
 (r_fd)(a,b)=\bigwedge_{\substack{fx=a\\fy=b}}d(x,y).
\]
These operations make $V_L$ noetherian. To see the first identity in \textup{(N1)}, a pseudometric zero on every fibre of $f$ is constant on pairs of fibres and hence descends to $f(G)$. The second identity follows by extending a pseudometric on $f(G)$ with distance $\infty$ between distinct points whenever at least one is outside $f(G)$. Kernels and images are the $0,\infty$-valued pseudometrics whose zero relations are, respectively,
\[
 R_f=\{(x,y):fx=fy\},\qquad
 \Delta_H\cup\bigl(f(G)\times f(G)\bigr).
\]
The usual quotient--subgroup factorization gives \textup{(N2)}. Joins of group congruences and intersections of subgroups give \textup{(N3)}. The same argument works over pointed sets, using pointed subsets in place of subgroups. None of these arguments requires distributivity of coefficient addition.

Change of coefficients also extends: the maps in $\MCL$ act by $d\mapsto h\circ d$, giving a faithful functor
\[
 Q:\MCL\longrightarrow\NF^c(\Grp),\qquad Q(L)=V_L.
\]
Indeed, the left adjoint of $h$ preserves $0$ and is subadditive, so the proof of Proposition~\ref{prop:coefficient} applies to pseudometrics. Faithfulness follows from the fibre over $C_2$, which is $L^{\op}$. In the idempotent case, sublevel equivalence relations give the analogue of \eqref{eq:tensor}:
\[
 V_L^{-1}(G)\cong L^{\op}\otimes_{\Sup}\operatorname{Eq}(|G|).
\]
These observations suggest studying pseudometrics as a common extension of the norm and equivalence-relation constructions.

Several stronger conclusions do not extend. For $C_4\to C_2$, the kernel has zero classes $\{0,2\}$ and $\{1,3\}$, whereas a conormal cluster has at most one nonsingleton zero class. Thus normal clusters need not be conormal. Moreover, $Q$ is not full: inversion induces an automorphism of forms
\[
 J_G(d)(x,y)=d(-x,-y).
\]
It commutes with inverse images, and hence with direct images, and preserves all fibre operations. It fixes the fibre over $C_2$ pointwise, but on $C_3$ it exchanges the $0,\infty$-valued pseudometrics with zero partitions $\{0,1\}\mid\{2\}$ and $\{0,2\}\mid\{1\}$. Consequently it cannot be induced by a coefficient map. Full faithfulness may therefore require additional structure or a different choice of base category and morphisms.

Even for quantale coefficients, the closure in direct images generally remains necessary. For example, on $C_6$ take the $0,\infty$-valued pseudometric with zero classes $\{0,1\}$, $\{4,5\}$, $\{2\}$ and $\{3\}$. Along reduction modulo $3$, its raw direct image has distances $0$ between $0,1$ and between $1,2$, but distance $\infty$ between $0,2$. It fails the triangle inequality. Thus Proposition~\ref{prop:quantaledirect} depends on the group-norm setting.

\medskip
\noindent\textbf{Adjunctions.}
Passing to arbitrary pseudometrics does not give $Q$ an adjoint in either direction with the morphisms of Definition~\ref{def:nfmorphism}. In fact, there are no morphisms of noetherian forms in either direction between $\mathsf S$ and any $V_L$. A morphism $\mathsf S\to V_L$ would preserve the equality between the image of $C_2=\{0,2\}\hookrightarrow C_4$ and the kernel of $C_4\to C_2$. In $V_L$ their zero partitions are, respectively,
\[
 \{0,2\}\mid\{1\}\mid\{3\},\qquad
 \{0,2\}\mid\{1,3\},
\]
which are different. Conversely, the fibre of $V_L$ over $C_3$ contains the bounded sublattice $M_3$ of $0,\infty$-valued pseudometrics, while $\mathsf S(C_3)=\two$. No bounded lattice homomorphism $M_3\to\two$ exists: among the three atoms, two would have the same value, contradicting preservation of their meet or join. An adjunction $A\dashv Q$ would supply the impossible unit $\mathsf S\to Q(A\mathsf S)$, and $Q\dashv B$ would supply the impossible counit $Q(B\mathsf S)\to\mathsf S$. This obstruction persists even if one restricts to forms with singleton fibre over the trivial group. In particular, the forms $V_L$ are not optimal, by Section~\ref{sec:optimal}.

There is, however, a reflection for each fixed coefficient lattice:
\[
 R_L:\mathbf{PMetGrp}_L\rightleftarrows\NGrp_L:I_L,
 \qquad R_L\dashv I_L,
\]
where $I_L(G,\nu)=(G,d_\nu)$ is fully faithful. Using the norm closure from \eqref{eq:closure}, put
\[
 u_d(g)=\bigwedge_{x\in G}d(x,x+g),\qquad
 R_L(G,d)=(G,c_G(u_d)).
\]
The corresponding invariant pseudometric is the greatest one pointwise below $d$. For a homomorphism $f:G\to H$ and an $L$-norm $\mu$ on $H$, the condition that $f:(G,d)\to I_L(H,\mu)$ be nonexpansive is equivalent to $\mu\circ f\leq u_d$, and hence to $\mu\circ f\leq c_G(u_d)$, since $\mu\circ f$ is a norm. This proves the adjunction. It is an adjunction of the total categories over $\Grp$, not an adjunction between coefficient lattices and noetherian forms. In particular, $I_L$ is not a morphism of noetherian forms: for the above inclusion $C_2\hookrightarrow C_4$, it sends the norm-form image to the coset pseudometric, rather than the image in $V_L$. A systematic account of this reflection, and of coefficient reconstruction under suitable additional requirements, is left for future work.

\enlargethispage{2\baselineskip}
\end{document}